\documentclass[11pt, twoside]{article}
\usepackage{latexsym}
\usepackage{amsmath}
\usepackage{amssymb}
\usepackage[all]{xy}
\usepackage{amsfonts}
\usepackage{verbatim}
\usepackage{amsthm}
\usepackage{bm}
\usepackage{mathrsfs}
\usepackage{epsfig}
\usepackage{xy}
\usepackage{array}
\usepackage{stmaryrd}
\usepackage{graphicx,color}
\usepackage{xcolor}
\usepackage{tikz}
\usetikzlibrary{arrows,calc}
\usepackage{etex}
\usepackage{mathdots}
\usepackage{float}
\usepackage{graphics}
\usepackage{pdflscape}
\usepackage{CJK}
\usepackage{anysize,hyperref}
\input xypic
\xyoption{all}
\usepackage{extarrows}
\usepackage[perpage,symbol]{footmisc}
\usepackage{setspace}
\def\C{\mathcal{C}}

\def\dr{\ar@{->}[r]}

\def\Hom{\mbox{Hom}}

\def\rad{\mbox{rad}}
\def\End{\mbox{End}}

\newcommand{\ind}{\mathsf{ind}\hspace{.01in}}
\newcommand{\supp}{\mathsf{Supp}\hspace{.01in}}
\begin{document}
\baselineskip=15pt
\title{\Large{\bf On the existence of Auslander-Reiten $\bm{(d+2)}$-angles in\\[2mm] $\bm{(d+2)}$-angulated categories$^\bigstar$\footnotetext{\hspace{-1em}$^\bigstar$Panyue Zhou was supported by the National Natural Science Foundation of China (Grant No. 11901190 and 11671221), and by the Hunan Provincial Natural Science Foundation of China (Grant No. 2018JJ3205), and by the Scientific Research Fund of Hunan Provincial Education Department (Grant No. 19B239).}}}
\medskip
\author{Panyue Zhou}

\date{}

\maketitle
\def\blue{\color{blue}}
\def\red{\color{red}}

\newtheorem{theorem}{Theorem}[section]
\newtheorem{lemma}[theorem]{Lemma}
\newtheorem{corollary}[theorem]{Corollary}
\newtheorem{proposition}[theorem]{Proposition}
\newtheorem{conjecture}{Conjecture}
\theoremstyle{definition}
\newtheorem{definition}[theorem]{Definition}
\newtheorem{question}[theorem]{Question}
\newtheorem{remark}[theorem]{Remark}
\newtheorem{remark*}[]{Remark}
\newtheorem{example}[theorem]{Example}
\newtheorem{example*}[]{Example}
\newtheorem{condition}[theorem]{Condition}
\newtheorem{condition*}[]{Condition}
\newtheorem{construction}[theorem]{Construction}
\newtheorem{construction*}[]{Construction}

\newtheorem{assumption}[theorem]{Assumption}
\newtheorem{assumption*}[]{Assumption}

\baselineskip=17pt
\parindent=0.5cm

\begin{abstract}
\baselineskip=16pt
Let $\C$ be a $(d+2)$-angulated category. In this note, we show that if $\C$ is a locally finite, then $\C$ has Auslander-Reiten $(d+2)$-angles. This extends a result of Xiao-Zhu for triangulated
categories.\\[0.5cm]
\textbf{Key words:} $(d+2)$-angulated categories; Auslander-Reiten $(d+2)$-angles; locally finite.\\[0.2cm]
\textbf{ 2010 Mathematics Subject Classification:} 16G70; 18E30.
\medskip
\end{abstract}

\pagestyle{myheadings}
\markboth{\rightline {\scriptsize   Panyue Zhou}}
         {\leftline{\scriptsize  On the existence of Auslander-Reiten $(d+2)$-angles in $(d+2)$-angulated categories}}

\section{Introduction}
Auslander-Reiten theory was introduced by Auslander and Reiten in \cite{ar1,ar2}. Since its introduction, Auslander-Reiten theory has become a fundamental tool for
studying the representation theory of Artin algebras.
Later it has been generalized to these situation of exact categories \cite{ji}, triangulated categories \cite{h,rv} and its subcategories \cite{as,j}
and some certain additive categories \cite{l,j,s} by many authors.
Extriangulated categories were recently introduced by Nakaoka and Palu \cite{np} as a simultaneous
generalization of exact categories and triangulated categories.
Hence, many results hold on exact categories and triangulated categories can be unified in the same framework.
Iyama, Nakaoka and Palu \cite{inp}  introduced the notion of almost split extensions and Auslander-Reiten-Serre duality for extriangulated categories, and gave explicit connections between these notions and also with the
classical notion of dualizing $k$-varieties.
Xiao and Zhu \cite{xz1,xz2} showed that if a triangulated category $\C$ is locally finite, then $\C$ has 
 Auslander-Reiten triangles. Recently, Zhu and Zhuang \cite{zz} proved that if an extriangulated category $\C$ is locally finite, then $\C$ has Auslander-Reiten $\mathbb{E}$-triangles.

In \cite{gko}, Geiss, Keller and Oppermann introduced $(d+2)$-angulated categories. These are
generalizations of triangulated categories, in the sense that triangles are replaced by
$(d+2)$-angles, that is, morphism sequences of length $d+2$. Thus a $1$-angulated category is precisely
a triangulated category.
Iyama and Yoshino \cite{iy} defined Auslander-Reiten $(d+2)$-angle in special $(d+2)$-angulated categories. Later,
Fedele \cite{f} defined Auslander-Reiten $(d+2)$-angles in additive subcategories of $(d+2)$-angulated categories closed
under $d$-extensions, an example of which are wide subcategories. He also proved that there are Auslander-Reiten $(d+2)$-angles in certain additive subcategories of $(d+2)$-angulated categories.
Recently, the author \cite{z} showed that a $(d+2)$-angulated category $\C$ has Auslander-Reiten $(d+2)$-angles if and only if
$\C$ has a Serre functor.

In this note, we continue to study Auslander-Reiten $(d+2)$-angles in $(d+2)$-angulated categories.
We will generalise Xiao-Zhu's result into $(d+2)$-angulated categories. Moreover, our proof is not far from the usual triangulated case.

Our main result is the following.
\begin{theorem}\emph{(see Theorem \ref{main} for details)}
Let $\C$ be a locally finite $(d+2)$-angulated category. If $X$ is an indecomposable, then there exists an Auslander-Reiten $(d+2)$-angle ending at $X$, and if $X$ is an indecomposable, then there exists an Auslander-Reiten $(d+2)$-angle starting at $X$. Thus $\C$ has Auslander-Reiten $(d+2)$-angles.
\end{theorem}

This article is organised as follows: In Section 2, we review some elementary definitions
that we need to use, including  $(d+2)$-angulated categories and Auslander-Reiten $(d+2)$ angles.
In Section 3,  we show our main result.

\section{Preliminaries}
In this section, we first recall the definition and basic properties of $(d+2)$-angulated categories from \cite{gko}.
Let $\C$ be an additive category with an automorphism $\Sigma^d:\C\rightarrow\C$, and an integer $d$ greater than or equal to one.

A $(d+2)$-$\Sigma^d$-$sequence$ in $\C$ is a sequence of objects and morphisms
$$A_0\xrightarrow{f_0}A_1\xrightarrow{f_1}A_2\xrightarrow{f_2}\cdots\xrightarrow{f_{d-1}}A_n\xrightarrow{f_d}A_{d+1}\xrightarrow{f_{d+1}}\Sigma^d A_0.$$
Its {\em left rotation} is the $(d+2)$-$\Sigma^d$-sequence
$$A_1\xrightarrow{f_1}A_2\xrightarrow{f_2}A_3\xrightarrow{f_3}\cdots\xrightarrow{f_{d}}A_{d+1}\xrightarrow{f_{d+1}}\Sigma^d A_0\xrightarrow{(-1)^{d}\Sigma^d f_0}\Sigma^d A_1.$$
A \emph{morphism} of $(d+2)$-$\Sigma^d$-sequences is  a sequence of morphisms $\varphi=(\varphi_0,\varphi_1,\cdots,\varphi_{d+1})$ such that the following diagram commutes
$$\xymatrix{
A_0 \ar[r]^{f_0}\ar[d]^{\varphi_0} & A_1 \ar[r]^{f_1}\ar[d]^{\varphi_1} & A_2 \ar[r]^{f_2}\ar[d]^{\varphi_2} & \cdots \ar[r]^{f_{d}}& A_{d+1} \ar[r]^{f_{d+1}}\ar[d]^{\varphi_{d+1}} & \Sigma^d A_0 \ar[d]^{\Sigma^d \varphi_0}\\
B_0 \ar[r]^{g_0} & B_1 \ar[r]^{g_1} & B_2 \ar[r]^{g_2} & \cdots \ar[r]^{g_{d}}& B_{d+1} \ar[r]^{g_{d+1}}& \Sigma^d B_0
}$$
where each row is a $(d+2)$-$\Sigma^d$-sequence. It is an {\em isomorphism} if $\varphi_0, \varphi_1, \varphi_2, \cdots, \varphi_{d+1}$ are all isomorphisms in $\C$.

\begin{definition}\cite[Definition 2.1]{gko}
A $(d+2)$-\emph{angulated category} is a triple $(\C, \Sigma^d, \Theta)$, where $\C$ is an additive category, $\Sigma^d$ is an automorphism of $\C$ ($\Sigma^d$ is called the $d$-suspension functor), and $\Theta$ is a class of $(d+2)$-$\Sigma^d$-sequences (whose elements are called $(d+2)$-angles), which satisfies the following axioms:
\begin{itemize}
\item[\textbf{(N1)}]
\begin{itemize}
\item[(a)] The class $\Theta$ is closed under isomorphisms, direct sums and direct summands.

\item[(b)] For each object $A\in\C$ the trivial sequence
$$ A\xrightarrow{1_A}A\rightarrow 0\rightarrow0\rightarrow\cdots\rightarrow 0\rightarrow \Sigma^dA$$
belongs to $\Theta$.

\item[(c)] Each morphism $f_0\colon A_0\rightarrow A_1$ in $\C$ can be extended to $(d+2)$-$\Sigma^d$-sequence: $$A_0\xrightarrow{f_0}A_1\xrightarrow{f_1}A_2\xrightarrow{f_2}\cdots\xrightarrow{f_{d-1}}A_d\xrightarrow{f_d}A_{d+1}\xrightarrow{f_{d+1}}\Sigma^d A_0.$$
\end{itemize}
\item[\textbf{(N2)}] A $(d+2)$-$\Sigma^n$-sequence belongs to $\Theta$ if and only if its left rotation belongs to $\Theta$.

\item[\textbf{(N3)}] Each solid commutative diagram
$$\xymatrix{
A_0 \ar[r]^{f_0}\ar[d]^{\varphi_0} & A_1 \ar[r]^{f_1}\ar[d]^{\varphi_1} & A_2 \ar[r]^{f_2}\ar@{-->}[d]^{\varphi_2} & \cdots \ar[r]^{f_{d}}& A_{d+1} \ar[r]^{f_{d+1}}\ar@{-->}[d]^{\varphi_{d+1}} & \Sigma^d A_0 \ar[d]^{\Sigma^d \varphi_0}\\
B_0 \ar[r]^{g_0} & B_1 \ar[r]^{g_1} & B_2 \ar[r]^{g_2} & \cdots \ar[r]^{g_{d}}& B_{d+1} \ar[r]^{g_{d+1}}& \Sigma^d B_0
}$$ with rows in $\Theta$, the dotted morphisms exist and give a morphism of  $(d+2)$-$\Sigma^d$-sequences.

\item[\textbf{(N4)}] In the situation of (N3), the morphisms $\varphi_2,\varphi_3,\cdots,\varphi_{d+1}$ can be chosen such that the mapping cone
$$A_1\oplus B_0\xrightarrow{\left(\begin{smallmatrix}
                                        -f_1&0\\
                                        \varphi_1&g_0
                                       \end{smallmatrix}
                                     \right)}
A_2\oplus B_1\xrightarrow{\left(\begin{smallmatrix}
                                        -f_2&0\\
                                        \varphi_2&g_1
                                       \end{smallmatrix}
                                     \right)}\cdots\xrightarrow{\left(\begin{smallmatrix}
                                        -f_{d+1}&0\\
                                        \varphi_{d+1}&g_d
                                       \end{smallmatrix}
                                     \right)} \Sigma^n A_0\oplus B_{d+1}\xrightarrow{\left(\begin{smallmatrix}
                                        -\Sigma^d f_0&0\\
                                        \Sigma^d\varphi_1&g_{d+1}
                                       \end{smallmatrix}
                                     \right)}\Sigma^dA_1\oplus\Sigma^d B_0$$
belongs to $\Theta$.
   \end{itemize}
\end{definition}

Now we give an example of $(d+2)$-angulated categories.

\begin{example}\label{def}
We recall the standard construction of $(d+2)$-angulated categories given by Geiss-Keller-Oppermann \cite[Theorem 1]{gko}.
Let $\C$ be a triangulated category and $\mathcal{T}$ a $d$-cluster tilting subcategory which is closed under $\Sigma^{d}$, where $\Sigma$ is the shift functor of $\C$. Then $(\mathcal{T},\Sigma^{d},\Theta)$ is a $(d+2)$-angulated category, where $\Theta$ is the class of all sequences
$$A_0\xrightarrow{f_0}A_1\xrightarrow{f_1}A_2\xrightarrow{f_2}\cdots\xrightarrow{f_{d-1}}A_d\xrightarrow{f_d}A_{d+1}\xrightarrow{f_{d+1}}\Sigma^{d} A_0$$
such that there exists a diagram
$$\xymatrixcolsep{0.3pc}
 \xymatrix{& A_1 \ar[dr]\ar[rr]^{f_1}  &  & A_2  \ar[dr]  & & \cdots  & & A_{d} \ar[dr]^{f_{d}}      \\
A_0 \ar[ur]^{f_0} & \mid & \ar[ll]  A_{1.5}\ar[ur] & \mid &  \ar[ll]  A_{2.5} & \cdots & A_{d-1.5}\ar[ur] & \mid & \ar[ll] A_{d+1}   }$$
with $A_i\in\mathcal{T}$ for all $i\in\mathbb{Z}$, such that all oriented triangles are triangles in $\C$, all non-oriented triangles commute, and $f_{d+1}$ is the composition along the lower edge of the diagram.
\end{example}

The following two lemmas are very useful which are needed later on.

\begin{lemma}\emph{\cite[Lemma 3.13]{f}}\label{y1}
Let $\C$ be a $(d+2)$-angulated category, and
\begin{equation}\label{t1}
\begin{array}{l}
A_0\xrightarrow{\alpha_0}A_1\xrightarrow{\alpha_1}A_2\xrightarrow{\alpha_2}\cdots\xrightarrow{\alpha_{d-1}}A_d\xrightarrow{\alpha_d}A_{d+1}\xrightarrow{\alpha_{d+1}}\Sigma^d A_0.\end{array}
\end{equation}
a $(d+2)$-angle in $\C$. Then the following are equivalent:
\begin{itemize}
\item[\rm (1)] $\alpha_0$ is a section;

\item[\rm (2)] $\alpha_d$ is a retraction;

\item[\rm (3)] $\alpha_{d+1}=0$.
\end{itemize}
If a $(d+2)$-angle \emph{(\ref{t1})} satisfies one of the above equivalent conditions, it is called \emph{split}.
\end{lemma}

\begin{lemma}\label{y4}\emph{\cite[Corollary 3.4]{lz}}
Let $\C$ be a $(d+2)$-angulated category,
and
$$A_0\xrightarrow{\alpha_0}A_1\xrightarrow{\alpha_1}A_2\xrightarrow{\alpha_2}\cdots\xrightarrow{\alpha_{d-1}}A_d\xrightarrow{\alpha_d}A_{d+1}\xrightarrow{\alpha_{d+1}}\Sigma^d A_0.$$
a $(d+2)$-angle in $\C$. Then for any morphism $\varphi_0\colon A_0\to B_0$, there exists
the following commutative diagram of $(d+2)$-angles
$$\xymatrix{
A_0 \ar[r]^{\alpha_0}\ar[d]^{\varphi_0}& A_1 \ar[r]^{\alpha_1}\ar@{-->}[d]^{\varphi_1} & A_2 \ar[r]^{\alpha_2}\ar@{-->}[d]^{\varphi_2}  & \cdots \ar[r]^{\alpha_{d-1}}& A_d \ar[r]^{\alpha_d}\ar@{-->}[d]^{\varphi_d}&A_{d+1}\ar[r]^{\alpha_{d+1}}\ar@{=}[d] & \Sigma^d A_0 \ar[d]^{\Sigma^d\varphi_0}\\
B_0 \ar[r]^{\beta_0}&B_1 \ar[r]^{\beta_1} & B_2 \ar[r]^{\beta_2}  & \cdots \ar[r]^{\beta_{d-1}} & B_d \ar[r]^{\beta_d}&A_{d+1}\ar[r]^{\beta_{d+1}}& \Sigma^d B_0\\
}$$
such that
$$A_0\xrightarrow{\left(
                    \begin{smallmatrix}
                      -\alpha_0 \\
                      \varphi_0 \\
                    \end{smallmatrix}
                  \right)} A_1\oplus B_0\xrightarrow{\left(
                             \begin{smallmatrix}
                               -\alpha_1 & 0 \\
                               \varphi_1 & \beta_0
                             \end{smallmatrix}
                           \right)}
 A_2\oplus B_1\xrightarrow{\left(
                                       \begin{smallmatrix}
                                         -\alpha_2 & 0  \\
                                         \varphi_2 &\beta_1
                                       \end{smallmatrix}
                                     \right)}\cdots\xrightarrow{\left(
                                       \begin{smallmatrix}
                                         -\alpha_{d-1}& 0  \\
                                         \varphi_{d-1} &\beta_{d-2}
                                       \end{smallmatrix}
                                     \right)}A_d\oplus B_{d-1}$$
                                     $\xrightarrow{\left(
                                       \begin{smallmatrix}
                                         \varphi_d,&\beta_{d-1}
                                       \end{smallmatrix}
                                     \right)}B_d\xrightarrow{(-1)^d\alpha_{d+1}\beta_d}\Sigma^dA_0$
 is a $(d+2)$-angle in $\C$.
\end{lemma}

Now we recall an Auslander-Reiten $(d+2)$ theory in $(d+2)$-angulated categories.

\setcounter{equation}{0}
We denote by ${\rm rad}_{\C}$ the Jacobson radical of $\C$. Namely, ${\rm rad}_{\C}$ is an ideal of $\C$ such that ${\rm rad}_{\C}(A, A)$
coincides with the Jacobson radical of the endomorphism ring ${\rm End}(A)$ for
any $A\in\C$.

\begin{definition}\cite[Definition 3.8]{iy} and \cite[Definition 5.1]{f}
Let $\C$ be a $(d+2)$-angulated category.
A $(d+2)$-angle
$$A_{\bullet}:~~\xymatrix {A_0 \xrightarrow{~\alpha_0~}A_1 \xrightarrow{~\alpha_1~} A_2 \xrightarrow{~\alpha_2~} \cdots
  \xrightarrow{~\alpha_{d - 1}~} A_d \xrightarrow{~\alpha_{d}~} A_{d+1}\xrightarrow{~\alpha_{d+1}~} \Sigma^d A_0}$$
in $\C$ is called an \emph{Auslander-Reiten $(d+2)$-angle }if
$\alpha_0$ is left almost split, $\alpha_d$ is right almost split and
when $d\geq 2$, also $\alpha_1,\alpha_2,\cdots,\alpha_{d-1}$ are in $\rad_{\C}$.
\end{definition}

\begin{remark}\cite[Remark 5.2]{f}
Assume $A_{\bullet}$ as in the above definition is an Auslander-Reiten $(d+2)$-angle.
Since $\alpha_0$ is left almost split implies that $\End(A_0)$ is local and hence $A_0$ is
indecomposable. Similarly, since $\alpha_d$ is right almost split, then $\End(A_{d+1})$ is local and hence $A_{d+1}$ is indecomposable.
Moreover, when $d=1$, we have $\alpha_0$ and $\alpha_d$ in $\rad_{\C}$, so that $\alpha_d$ is right minimal and $\alpha_0$ is left minimal. When $d\geq 2$, since $\alpha_{d-1}\in\rad_{\C}$, we have that $\alpha_d$ is right minimal and similarly $\alpha_0$ is left
minimal.
\end{remark}

\begin{remark}\cite[Lemma 5.3]{f}
Let $\C$ be a $(d+2)$-angulated category and
$$A_{\bullet}:~~\xymatrix {A_0 \xrightarrow{~\alpha_0~}A_1 \xrightarrow{~\alpha_1~} A_2 \xrightarrow{~\alpha_2~} \cdots
  \xrightarrow{~\alpha_{d - 1}~} A_d \xrightarrow{~\alpha_{d}~} A_{d+1}\xrightarrow{~\alpha_{d+1}~} \Sigma^d A_0}$$
be a $(d+2)$-angle in $\C$. Then the following statements are equivalent:
\begin{itemize}
\item[(1)] $A_{\bullet}$ is an Auslander-Reiten $(d+2)$-angle;
\item[(2)] $\alpha_0,\alpha_1,\cdots,\alpha_{d-1}$ are in ${\rm rad}_{\C}$ and $\alpha_d$ is right almost split;
\item[(3)] $\alpha_1,\alpha_2,\cdots,\alpha_{d}$ are in ${\rm rad}_{\C}$ and $\alpha_0$ is left almost split.
\end{itemize}
\end{remark}

\begin{lemma}\emph{\cite[Lemma 5.4]{f}}\label{lem0}
Let $\C$ be a $(d+2)$-angulated category and
$$A_{\bullet}:~~\xymatrix {A_0 \xrightarrow{~\alpha_0~}A_1 \xrightarrow{~\alpha_1~} A_2 \xrightarrow{~\alpha_2~} \cdots
  \xrightarrow{~\alpha_{d - 1}~} A_d \xrightarrow{~\alpha_{d}~} A_{d+1}\xrightarrow{~\alpha_{d+1}~} \Sigma^d A_0}$$
be a $(d+2)$-angle in $\C$.
Assume that $\alpha_d$ is right almost split and if $d\geq2$, also that $\alpha_1,\alpha_2,\cdots,\alpha_{d-1}$ are in $\emph{\rad}_{\C}$.
Then the following are equivalent:
\begin{itemize}
\item[\emph{(1)}] $A_{\bullet}$ is an Auslander-Reiten $(d+2)$-angle;
\item[\emph{(2)}] $\emph{\End}(A_0)$ is local;
\item[\emph{(3)}] $\alpha_{d+1}$ is left minimal;
\item[\emph{(4)}] $\alpha_{0}$ is in $\emph{\rad}_{\C}$.
\end{itemize}
\end{lemma}

In the case $d=1$, so in the case of a triangulated category, a morphism can be extended to a triangle
in a unique way up to isomorphism. On the other hand, for $d\geq 2$, a morphism can be extendend to a
$(d+2)$-angle in different non-isomorphic ways. However, we still have a unique ``minimal" $(d+2)$-angle
extending any given morphism.

\begin{lemma}\emph{\cite[Lemma 5.18]{ot} and \cite[Lemma 3.14]{f}}\label{y3}
Let $d\geq 2$ and $h\colon A_{d+1}\to \Sigma^d A_0$ be any morphism in a $(d+2)$-angulated category $\C$. Then, up to
isomorphism, there exists a unique $(d+2)$-angle of the form
$$\xymatrix {A_0 \xrightarrow{~\alpha_0~}A_1 \xrightarrow{~\alpha_1~} A_2 \xrightarrow{~\alpha_2~} \cdots
  \xrightarrow{~\alpha_{d-2}~} A_{d-1} \xrightarrow{~\alpha_{d - 1}~} A_d \xrightarrow{~\alpha_{d}~} A_{d+1}\xrightarrow{~h~} \Sigma^d A_0}$$
with $\alpha_1,\alpha_2,\cdots,\alpha_{d-1}$ in ${\rm rad_{\C}}$.
\end{lemma}

\section{Proof of main result}
In this section, let $k$ be a field. We always assume that $\C$ is a $k$-linear Hom-finite Krull-Schmidt $(d+2)$-angulated category.
We denote by $\ind(\C)$ the set of isomorphism classes of indecomposable objects in $\C$.
For any $X\in\ind(\C)$, we denote by $\supp \Hom_{\C}(X,-)$ the subcategory of $\C$ generated by
objects $Y$ in $\ind(\C)$ with $\Hom_{\C}(X,Y)\neq0$. Similarly, $\supp\Hom_{\C}(-,X)$ denotes
the subcategory generated by objects $Y$ in $\ind(\C)$ with $\Hom_{\C}(Y,X)\neq 0$. If
$\supp\Hom_{\C}(X,-)$ ($\supp\Hom_{\C}(-,X)$, respectively) contains only finitely many
indecomposables, we say that $|\supp\Hom_{\C}(X,-)|<\infty $ ($|\supp\Hom_{\C}(-,X)|<\infty$
respectively).
\medskip

Based on the definition of locally finite triangulated categories \cite{xz1,xz2},
we define the notion of locally finite $(d+2)$-angulated categories.
\begin{definition}
A $(d+2)$-angulated category $\C$ is called \emph{locally finite}
if $|\supp\Hom_{\C}(X,-)|<\infty$ and $|\supp\Hom_{\C}(-,X)|<\infty$, for any object $X\in\ind(C)$.
\end{definition}

We know that the derived categories of finite dimensional hereditary algebras of finite type and the
stable module categories of finite dimensional self-injective algebras of finite type are examples of locally finite triangulated categories, see \cite{xz1,xz2}. In those locally finite triangulated categories, we take a $d$-cluster titling subcategory which is closed under the $d$-th power of the shift functor. By Example \ref{def},  we can obtain locally finite $(d+2)$-angulated categories.

\begin{definition}\label{y2}
Let $\C$ be a $(d+2)$-angulated category and $X\in\ind(\C)$. We define a set of $(d+2)$-angles as follows:
$$
S(X):=\biggl\{A_{\bullet}:A\xrightarrow{\alpha_0}A_1\xrightarrow{\alpha_1}\cdots\xrightarrow{\alpha_{d-1}}A_{d}\xrightarrow{\alpha_d}X\xrightarrow{\alpha_{d+1}}\Sigma^d A_0~\bigg|~\begin{array}{ll}
 A_{\bullet}~\mbox{is a non-split $(d+2)$-angle}\\
\textrm{with}~A\in\ind(\C), \mbox{and when}\\ d\geq2,~\alpha_1,\alpha_2,\cdots,\alpha_{d-1}$ in ${\rm rad_{\C}}.
 \end{array}\biggl\}$$

Dually, we can define a set of $(d+2)$-angles as follows:
$$
T(X):=\biggl\{A_{\bullet}:X\xrightarrow{\alpha_0}A_1\xrightarrow{\alpha_1}\cdots\xrightarrow{\alpha_{d-1}}A_{d}\xrightarrow{\alpha_d}A\xrightarrow{\alpha_{d+1}}\Sigma^d A_0~\bigg|~\begin{array}{ll}
 A_{\bullet}~\mbox{is a non-split $(d+2)$-angle}\\
\textrm{with}~A\in\ind(\C), \mbox{and when}\\ d\geq2,~ \alpha_1,\alpha_2,\cdots,\alpha_{d-1}$ in ${\rm rad_{\C}}.
 \end{array}\biggl\}$$
\end{definition}

\begin{lemma}\label{y5}
Let $\C$ be a $(d+2)$-angulated category and $X\in\ind(\C)$. Then $S(X)$ and $T(X)$
are non-empty.\end{lemma}

\proof We only show that $S(X)$ non-empty, dually one can show that $T(X)$ is non-empty.

Since $X\in\ind(\C)$, there is an object $A\in\C$ such that $\Hom_{\C}(X,\Sigma^dA)\neq 0$.
Thus there exists a non-split $(d+2)$-angle:
$$\xymatrix {B_{\bullet}:~A \xrightarrow{~\alpha_0~}B_1 \xrightarrow{~\alpha_1~} B_2 \xrightarrow{~\alpha_2~} \cdots
  \xrightarrow{~\alpha_{d-2}~} B_{d-1} \xrightarrow{~\alpha_{d - 1}~} B_d \xrightarrow{~\alpha_{d}~} X\xrightarrow{~h~} \Sigma^d A}.$$
We decompose $A$ into a direct sum of indecomposable objects $A=\bigoplus\limits_{i=1}^nA_i$.
Without loss of generality, we can assume that $A=U\oplus V$ where $U$ and $V$ are indecomposable.
By Lemma \ref{y4}, for the morphism $(1,0)\colon U\oplus V\to U$, there exists
the following commutative diagram of $(d+2)$-angles
$$\xymatrix@C=1.2cm{
 U\oplus V \ar[r]^{\;\;(u,~v)}\ar[d]^{(1,~0)}& A_1 \ar[r]^{\alpha_1}\ar@{-->}[d]^{\varphi_1} & A_2 \ar[r]^{\alpha_2}\ar@{-->}[d] & \cdots \ar[r]^{\alpha_{d-1}}& A_d \ar[r]^{\alpha_d}\ar@{-->}[d]&X\ar[r]^{h\qquad}\ar@{=}[d] & \Sigma^d U\oplus \Sigma^d V \ar[d]^{(1,~0)}\\
U \ar[r]^{\beta_0}&C_1 \ar[r]^{\beta_1} & C_2 \ar[r]^{\beta_2}  & \cdots \ar[r]^{\beta_{d-1}} & C_d \ar[r]^{\beta_d}&X\ar[r]^{\beta_{d+1}}& \Sigma^d U.
}$$
Similarly, for the morphism $(0,1)\colon U\oplus V\to V$, there exists
the following commutative diagram of $(d+2)$-angles
$$\xymatrix@C=1.2cm{
 U\oplus V \ar[r]^{\;\;(u,~v)}\ar[d]^{(0,~1)}& A_1 \ar[r]^{\alpha_1}\ar@{-->}[d]^{\psi_1} & A_2 \ar[r]^{\alpha_2}\ar@{-->}[d] & \cdots \ar[r]^{\alpha_{d-1}}& A_d \ar[r]^{\alpha_d}\ar@{-->}[d]&X\ar[r]^{h\qquad}\ar@{=}[d] & \Sigma^d U\oplus \Sigma^d V \ar[d]^{(0,~1)}\\
V \ar[r]^{\gamma_0}&D_1 \ar[r]^{\gamma_1} & D_2 \ar[r]^{\gamma_2}  & \cdots \ar[r]^{\gamma_{d-1}} & D_d \ar[r]^{\gamma_d}&X\ar[r]^{\gamma_{d+1}}& \Sigma^d V.
}$$
We claim that the at least one of the following two $(d+2)$-angles is non-split
$$\xymatrix{
U \ar[r]^{\beta_0}&C_1 \ar[r]^{\beta_1} & C_2 \ar[r]^{\beta_2}  & \cdots \ar[r]^{\beta_{d-1}} & C_d \ar[r]^{\beta_d}&X\ar[r]^{\beta_{d+1}\;\;}& \Sigma^d U,}$$
$$\xymatrix{
V \ar[r]^{\gamma_0}&D_1 \ar[r]^{\gamma_1} & D_2 \ar[r]^{\gamma_2}  & \cdots \ar[r]^{\gamma_{d-1}} & D_d \ar[r]^{\gamma_d}&X\ar[r]^{\gamma_{d+1}\;\;}& \Sigma^d V.
}$$
If not like this, by Lemma \ref{y1}, we obtain that $\beta_{d+1}=0=\gamma_{d+1}$. By (N3), we have the following commutative diagram of $(d+2)$-angles:
$$\xymatrix@C=0.9cm{
 U\oplus V \ar[r]^{\;\;(u,~v)}\ar@{=}[d]& A_1 \ar[r]^{\alpha_1}\ar[d]^{\left(\begin{smallmatrix}
                               \varphi_1\\ \psi_1
                             \end{smallmatrix}
                           \right)} & A_2 \ar[r]^{\alpha_2}\ar@{-->}[d] & \cdots \ar[r]^{\alpha_{d-1}}& A_d \ar[r]^{\alpha_d}\ar@{-->}[d]&X\ar[r]^{h\qquad}\ar@{-->}[d] & \Sigma^d U\oplus \Sigma^d V \ar@{=}[d]\\
U\oplus V \ar[r]^{\delta_0}&C_1\oplus D_1\ar[r]^{\delta_1} & C_2\oplus D_2 \ar[r]^{\delta_2}  & \cdots \ar[r]^{\delta_{d-1}} & C_d\oplus D_d \ar[r]^{\delta_d}&X\oplus X\ar[r]^{\delta_{d+1}\;\;\;\;\;\;}& \Sigma^d U\oplus\Sigma^dV
}$$
where $\delta_i=\left(\begin{smallmatrix}
                               \beta_i&0\\ 0&\gamma_i
                             \end{smallmatrix}
                           \right)$.
It follows that $h=0$. This is a contradiction since $B_{\bullet}$ is non-split.

For the morphism $\beta_{d+1}\neq 0$ or $\gamma_{d+1}\neq 0$, by Lemma \ref{y3}, we can find a $(d+2)$-angle as we want. This shows that $S(X)$ is nonempty.  \qed

\begin{definition}\label{y6}
Let $\C$ be a $(d+2)$-angulated category, and
$$A_{\bullet}:~A\xrightarrow{\alpha_0}A_1\xrightarrow{\alpha_1}\cdots\xrightarrow{\alpha_{d-1}}A_{d}\xrightarrow{\alpha_d}X\xrightarrow{\alpha_{d+1}}\Sigma^d A_0$$
$$B_{\bullet}:~B\xrightarrow{\beta_0}B_1\xrightarrow{\beta_1}\cdots\xrightarrow{\beta_{d-1}}B_{d}\xrightarrow{\beta_d}X\xrightarrow{\beta_{d+1}}\Sigma^d B_0$$
two $(d+2)$-angles in $S(X)$. We say that
$A_{\bullet}>B_{\bullet}$ if there are morphisms $\varphi_i\in\Hom_{\C}(A_i,B_i)$, $(i=0,1,\cdots,d)$ such that
the following diagram commutative:
$$\xymatrix{
A\ar[r]^{\alpha_0}\ar@{-->}[d]^{\varphi_0}& A_1 \ar[r]^{\alpha_1}\ar@{-->}[d]^{\varphi_1} & A_2 \ar[r]^{\alpha_2}\ar@{-->}[d]^{\varphi_2}  & \cdots \ar[r]^{\alpha_{d-1}}& A_d \ar[r]^{\alpha_d\;\;}\ar@{-->}[d]^{\varphi_d}&X\ar[r]^{\alpha_{d+1}}\ar@{=}[d] & \Sigma^d A_0 \ar@{-->}[d]^{\Sigma^d\varphi_0}\\
B\ar[r]^{\beta_0}&B_1 \ar[r]^{\beta_1} & B_2 \ar[r]^{\beta_2}  & \cdots \ar[r]^{\beta_{d-1}} & B_d \ar[r]^{\beta_d\;\;}&X\ar[r]^{\beta_{d+1}}& \Sigma^d B_0.
}$$
We say that $A_{\bullet}\sim B_{\bullet}$ if $\varphi_0$ is an isomorphism.

 Dually, let
$$A_{\bullet}:~X\xrightarrow{\alpha_0}A_1\xrightarrow{\alpha_1}\cdots\xrightarrow{\alpha_{d-1}}A_{d}\xrightarrow{\alpha_d}A\xrightarrow{\alpha_{d+1}}\Sigma^d A_0$$
$$B_{\bullet}:~X\xrightarrow{\beta_0}B_1\xrightarrow{\beta_1}\cdots\xrightarrow{\beta_{d-1}}B_{d}\xrightarrow{\beta_d}B\xrightarrow{\beta_{d+1}}\Sigma^d B_0$$
be two $(d+2)$-angles in $T(X)$.
 We say that
$A_{\bullet}>B_{\bullet}$ if there are morphisms $\varphi_i\in\Hom_{\C}(A_i,B_i)$, $(i=1,2,\cdots,d+1)$ such that
the following diagram commutative:
$$\xymatrix{
X\ar[r]^{\alpha_0}\ar@{=}[d]& A_1 \ar[r]^{\alpha_1}\ar@{-->}[d]^{\varphi_1} & A_2 \ar[r]^{\alpha_2}\ar@{-->}[d]^{\varphi_2}  & \cdots \ar[r]^{\alpha_{d-1}}& A_d \ar[r]^{\alpha_d\;\;}\ar@{-->}[d]^{\varphi_d}&A\ar[r]^{\alpha_{d+1}}\ar@{-->}[d]^{\varphi_{d+1}} & \Sigma^dX \ar@{=}[d]\\
 X\ar[r]^{\beta_0}&B_1 \ar[r]^{\beta_1} & B_2 \ar[r]^{\beta_2}  & \cdots \ar[r]^{\beta_{d-1}} & B_d \ar[r]^{\beta_d\;\;}&B\ar[r]^{\beta_{d+1}}& \Sigma^d X.
}$$
We say that $A_{\bullet}\sim B_{\bullet}$ if $\varphi_{d+1}$ is an isomorphism.
\end{definition}

\begin{lemma}\label{y7}
$S(X)$ is a direct ordered set with the relation defined in Definition \ref{y6}, and
$T(X)$ is a direct ordered set with the relation defined in Definition \ref{y6}.
\end{lemma}

\proof We just prove the first statement, the second statement proves similarly.

Assume that
$$A_{\bullet}:~A\xrightarrow{\alpha_0}A_1\xrightarrow{\alpha_1}\cdots\xrightarrow{\alpha_{d-1}}A_{d}\xrightarrow{\alpha_d}X\xrightarrow{\alpha_{d+1}}\Sigma^d A_0$$
$$B_{\bullet}:~B\xrightarrow{\beta_0}B_1\xrightarrow{\beta_1}\cdots\xrightarrow{\beta_{d-1}}B_{d}\xrightarrow{\beta_d}X\xrightarrow{\beta_{d+1}}\Sigma^d B_0$$
belong to $S(X)$.

We first show that if $A_{\bullet}>B_{\bullet}$ and $B_{\bullet}>A_{\bullet}$, then
$A_{\bullet}\sim B_{\bullet}$.

Since $A_{\bullet}>B_{\bullet}$ and $B_{\bullet}>A_{\bullet}$, we have the following two commutative diagrams
$$\xymatrix{
A\ar[r]^{\alpha_0}\ar@{-->}[d]^{\varphi_0}& A_1 \ar[r]^{\alpha_1}\ar@{-->}[d]^{\varphi_1} & A_2 \ar[r]^{\alpha_2}\ar@{-->}[d]^{\varphi_2}  & \cdots \ar[r]^{\alpha_{d-1}}& A_d \ar[r]^{\alpha_d\;\;}\ar@{-->}[d]^{\varphi_d}&X\ar[r]^{\alpha_{d+1}}\ar@{=}[d] & \Sigma^d A_0 \ar@{-->}[d]^{\Sigma^d\varphi_0}\\
B\ar[r]^{\beta_0}&B_1 \ar[r]^{\beta_1} & B_2 \ar[r]^{\beta_2}  & \cdots \ar[r]^{\beta_{d-1}} & B_d \ar[r]^{\beta_d\;\;}&X\ar[r]^{\beta_{d+1}}& \Sigma^d B_0,
}$$
$$\xymatrix{
B\ar[r]^{\beta_0}\ar@{-->}[d]^{\psi_0}& B_1 \ar[r]^{\beta_1}\ar@{-->}[d]^{\psi_1} & B_2 \ar[r]^{\beta_2}\ar@{-->}[d]^{\psi_2}  & \cdots \ar[r]^{\beta_{d-1}}& B_d \ar[r]^{\beta_d\;\;}\ar@{-->}[d]^{\psi_d}&X\ar[r]^{\beta_{d+1}}\ar@{=}[d] & \Sigma^d B_0 \ar@{-->}[d]^{\Sigma^d\psi_0}\\
A\ar[r]^{\alpha_0}&A_1 \ar[r]^{\alpha_1} & A_2 \ar[r]^{\alpha_2}  & \cdots \ar[r]^{\alpha_{d-1}} & A_d \ar[r]^{\alpha_d\;\;}&X\ar[r]^{\alpha_{d+1}}& \Sigma^d A_0.
}$$
Since $A$ is an indecomposable, we have that $\End(A)$ is local implies that
$\psi_0\varphi_0$ is nilpotent or is
an isomorphism. If $\psi_0\varphi_0$ is nilpotent, there exists a positive integer $m$ such that
$(\psi_0\varphi_0)^m=0$. We write $\omega_i=\psi_i\varphi_i$.
Thus we have the following commutative diagram
$$\xymatrix{
A\ar[r]^{\alpha_0}\ar[d]^{(\psi_0\varphi_0)^m}& A_1 \ar[r]^{\alpha_1}\ar[d]^{(\omega_1)^m} & A_2 \ar[r]^{\alpha_2}\ar[d]^{(\omega_2)^m}  & \cdots \ar[r]^{\alpha_{d-2}}&\ar[r]^{\alpha_{d-1}}\ar[d]^{(\omega_{d-1})^m}A_{d-1}&A_{d} \ar[r]^{\alpha_d\;\;}\ar[d]^{(\varphi_d)^m}&X\ar[r]^{\alpha_{d+1}}\ar@{=}[d] & \Sigma^d A_0 \ar[d]^{\Sigma^d(\psi_0\varphi_0)^m}\\
A\ar[r]^{\alpha_0}&A_1 \ar[r]^{\alpha_1} & A_2 \ar[r]^{\alpha_2}  & \cdots \ar[r]^{\alpha_{d-2}} &A_{d-1}\ar[r]^{\alpha_{d-1}}& A_d \ar[r]^{\alpha_d\;\;}&X\ar[r]^{\alpha_{d+1}}& \Sigma^d A_0.}$$
Then $\alpha_{d+1}=\Sigma^d(\psi_0\varphi_0)^m\alpha_{d+1}=0$. This is a contradiction since $A_{\bullet}$ is non-split.
Hence $\psi_0\varphi_0$ is an isomorphism.
By a similar argument we obtain that $\varphi_{0}\psi_{0}$ is an isomorphism.
This shows that $\varphi_0$ is isomorphism. So $A_{\bullet}\sim B_{\bullet}$.

It is clear that if $A_{\bullet}>B_{\bullet}$ and $B_{\bullet}>C_{\bullet}$, then
$A_{\bullet}\sim C_{\bullet}$.

Now we show that if $A_{\bullet},B_{\bullet}\in S(X)$, then there exists $C_{\bullet}\in S(X)$ such that $A_{\bullet}>C_{\bullet}$ and
$B_{\bullet}\sim C_{\bullet}$.

For the morphism $\beta_{d}\colon B_d\to X$, by the dual of Lemma \ref{y4},
there exists
the following commutative diagram of $(d+2)$-angles
$$\xymatrix@C=1.2cm{
A\ar[r]^{\gamma_0}\ar@{=}[d]& D_1 \ar[r]^{\gamma_1}\ar@{-->}[d]^{\psi_1} & D_2 \ar[r]^{\gamma_2}\ar@{-->}[d]^{\psi_2}  & \cdots \ar[r]^{\gamma_{d-1}}& D_d \ar[r]^{\gamma_d\;\;}\ar@{-->}[d]^{\psi_d}&B_d\ar[r]^{\gamma_{d+1}}\ar[d]^{\gamma_d} & \Sigma^d A \ar@{=}[d]\\
A\ar[r]^{\alpha_0}&A_1 \ar[r]^{\alpha_1} & A_2 \ar[r]^{\alpha_2}  & \cdots \ar[r]^{\alpha_{d-1}} & A_d \ar[r]^{\alpha_d}&X\ar[r]^{\alpha_{d+1}}& \Sigma^d A
}$$
such that
$$M_{\bullet}:~D_1\xrightarrow{} M_1\xrightarrow{}
 M_2\xrightarrow{}\cdots\xrightarrow{}M_{d-1}
                                     \xrightarrow{}B_d\oplus A_d\xrightarrow{\left(
                                       \begin{smallmatrix}
                                         \beta_d,&\alpha_{d}
                                       \end{smallmatrix}
                                     \right)}X\xrightarrow{~h~}\Sigma^dD_1$$
 is a $(d+2)$-angle in $\C$, where $M_i=D_{i+1}\oplus A_{i}~(i=1,2,\cdots,d-1)$.
Since $\beta_d$ and $\alpha_d$ are not retraction, we have that
$(\beta_d,\alpha_d)$ is also not retraction.
If not like this, there exists a morphism $\binom{s}{t}\colon X\to B_d\oplus A_d$
such that $(\beta_d,\alpha_d)\binom{s}{t}=1_{X}$ and then $\beta_ds+\alpha_dt=1_X$.
Since $X$ is an indecomposable, we have that $\End(X)$ is local implies that
either $\beta_ds$ or $\alpha_dt$ is an isomorphism.
Thus either $\beta_d$ or $\alpha_d$ is an retraction, a contradiction.
That is, $M_{\bullet}$ is non-split.

Without loss of generality, we can assume that $D_1=U\oplus V$ where $U$ and $V$ are indecomposable.
By Lemma \ref{y4}, for the morphism $(1,0)\colon U\oplus V\to U$, there exists
the following commutative diagram of $(d+2)$-angles
$$\xymatrix@C=1.2cm{
 U\oplus V \ar[r]^{\;\;(u,~v)}\ar[d]^{(1,~0)}& M_1 \ar[r]\ar@{-->}[d]^{\varphi_1} & M_2 \ar[r]\ar@{-->}[d] & \cdots \ar[r]&B_d\oplus A_d \ar[r]\ar@{-->}[d]&X\ar[r]\ar@{=}[d] & \Sigma^d U\oplus \Sigma^d V \ar[d]^{(1,~0)}\\
U \ar[r]^{\delta_0}&L_1 \ar[r]& L_2 \ar[r] & \cdots \ar[r] & L_d \ar[r]&X\ar[r]^{h}& \Sigma^d U.
}$$
Similarly, for the morphism $(0,1)\colon U\oplus V\to V$, there exists
the following commutative diagram of $(d+2)$-angles
$$\xymatrix@C=1.2cm{
 U\oplus V \ar[r]^{\;\;(u,~v)}\ar[d]^{(0,~1)}& M_1 \ar[r]\ar@{-->}[d]^{\psi_1} & M_2 \ar[r]\ar@{-->}[d] & \cdots \ar[r]&B_d\oplus A_d \ar[r]\ar@{-->}[d]&X\ar[r]\ar@{=}[d] & \Sigma^d U\oplus \Sigma^d V \ar[d]^{(0,~1)}\\
V \ar[r]^{\eta_0}&N_1 \ar[r]& N_2 \ar[r] & \cdots \ar[r] & N_d \ar[r]&X\ar[r]& \Sigma^d V.
}$$
Using similar arguments as in the proof of Lemma \ref{y5}, we conclude
that the at least one of the following two $(d+2)$-angles is non-split
$$\xymatrix{
U \ar[r]^{\delta_0}&L_1 \ar[r] & L_2 \ar[r]  & \cdots \ar[r] & L_d \ar[r]&X\ar[r]^{h\;\;\;}& \Sigma^d U,}$$
$$\xymatrix{
V \ar[r]^{\eta_0}&N_1 \ar[r] & N_2 \ar[r]  & \cdots \ar[r] & N_d \ar[r]&X\ar[r]& \Sigma^d V.
}$$
Without loss of generality, we assume that
$$\xymatrix{U \ar[r]^{\delta_0}&L_1 \ar[r] & L_2 \ar[r]  & \cdots \ar[r] & L_d \ar[r]&X\ar[r]^{h\;\;\;}& \Sigma^d U,}$$ is non-split.
By Lemma \ref{y3}, we can find a non-split $(d+2)$-angle
$$\xymatrix{
C_{\bullet}:~U \ar[r]^{\quad\lambda_0}&C_1 \ar[r]^{\lambda_1} & C_2 \ar[r]^{\lambda_2}  & \cdots \ar[r]^{\lambda_{d-1}} & C_d \ar[r]^{\lambda_d}&X\ar[r]^{h\;\;}& \Sigma^d U}$$
with $\lambda_1,\lambda_2,\cdots,\lambda_{d-1}$ in ${\rm rad_{\C}}$.
By (N3), we have the following commutative diagram
$$\xymatrix@C=1.2cm{A\ar[r]^{\alpha_0}\ar@{-->}[d]& A_1 \ar[r]^{\alpha_1}\ar@{-->}[d] & A_2 \ar[r]^{\alpha_2}\ar@{-->}[d]& \cdots \ar[r]^{\alpha_{d-1}}& A_d \ar[r]^{\alpha_d\;\;}\ar[d]^{\binom{0}{1}}&X\ar[r]^{\alpha_{d+1}}\ar@{=}[d] & \Sigma^d A_0 \ar@{-->}[d]\\
 U\oplus V \ar[r]^{\;\;(u,~v)}\ar[d]^{(1,~0)}& M_1 \ar[r]\ar[d]^{\varphi_1} & M_2 \ar[r]\ar[d] & \cdots \ar[r]&B_d\oplus A_d \ar[r]^{\quad(\beta_d,\alpha_d)}\ar[d]&X\ar[r]\ar@{=}[d] & \Sigma^d U\oplus \Sigma^d V \ar[d]^{(1,~0)}\\
U \ar[r]^{\delta_0}\ar@{=}[d] &L_1 \ar[r]\ar@{-->}[d] & L_2 \ar[r]\ar@{-->}[d]  & \cdots \ar[r] & L_d\ar@{-->}[d] \ar[r]&X\ar[r]^{h}\ar@{=}[d]& \Sigma^d U\ar@{=}[d]\\
U \ar[r]^{\quad\lambda_0}&C_1 \ar[r]^{\lambda_1} & C_2 \ar[r]^{\lambda_2}  & \cdots \ar[r]^{\lambda_{d-1}} & C_d \ar[r]^{\lambda_d}&X\ar[r]^{h\;\;}& \Sigma^d U
}$$
of $(d+2)$-angles. This shows that $A_{\bullet}>C_{\bullet}$.

By (N3), we have the following commutative diagram
$$\xymatrix@C=1.2cm{B\ar[r]^{\beta_0}\ar@{-->}[d]& B_1 \ar[r]^{\beta_1}\ar@{-->}[d] & B_2 \ar[r]^{\beta_2}\ar@{-->}[d]& \cdots \ar[r]^{\beta_{d-1}}& B_d \ar[r]^{\beta_d\;\;}\ar[d]^{\binom{1}{0}}&X\ar[r]^{\beta_{d+1}}\ar@{=}[d] & \Sigma^d B_0 \ar@{-->}[d]\\
 U\oplus V \ar[r]^{\;\;(u,~v)}\ar[d]^{(1,~0)}& M_1 \ar[r]\ar[d]^{\varphi_1} & M_2 \ar[r]\ar[d] & \cdots \ar[r]&B_d\oplus A_d \ar[r]^{\quad(\beta_d,\alpha_d)}\ar[d]&X\ar[r]\ar@{=}[d] & \Sigma^d U\oplus \Sigma^d V \ar[d]^{(1,~0)}\\
U \ar[r]^{\delta_0}\ar@{=}[d] &L_1 \ar[r]\ar@{-->}[d] & L_2 \ar[r]\ar@{-->}[d]  & \cdots \ar[r] & L_d\ar@{-->}[d] \ar[r]&X\ar[r]^{h}\ar@{=}[d]& \Sigma^d U\ar@{=}[d]\\
U \ar[r]^{\quad\lambda_0}&C_1 \ar[r]^{\lambda_1} & C_2 \ar[r]^{\lambda_2}  & \cdots \ar[r]^{\lambda_{d-1}} & C_d \ar[r]^{\lambda_d}&X\ar[r]^{h\;\;}& \Sigma^d U
}$$ of $(d+2)$-angles.
This shows that $B_{\bullet}>C_{\bullet}$. \qed
\medskip

\begin{lemma}\label{y8}
Let $\C$ be a locally finite $(d+2)$-angulated category and $X\in\ind(\C)$. Then
$S(X)$ has a minimal element, and
$T(X)$ has a minimal element.
\end{lemma}

\proof We just prove the first statement, the second statement proves similarly.

Since $X\in\ind(\C)$, there is an object $A\in\C$ such that $\Hom_{\C}(X,\Sigma^dA)\neq 0$.
Then there exists a non-split $(d+2)$-angle:
$$\xymatrix {A_{\bullet}:~A \xrightarrow{~\alpha_0~}A_1 \xrightarrow{~\alpha_1~} A_2 \xrightarrow{~\alpha_2~} \cdots
  \xrightarrow{~\alpha_{d-2}~} A_{d-1} \xrightarrow{~\alpha_{d - 1}~} A_d \xrightarrow{~\alpha_{d}~} X\xrightarrow{~h~} \Sigma^d A}.$$
We decompose $B_d$ into a direct sum of indecomposable objects $A_d=\bigoplus\limits_{k=1}^nB_k$.
Thus $A_{\bullet}$ can be written as
$$\xymatrix {A_{\bullet}:~A \xrightarrow{~\alpha_0~}A_1 \xrightarrow{~\alpha_1~} A_2 \xrightarrow{~\alpha_2~} \cdots
  \xrightarrow{~\alpha_{d-2}~} A_{d-1} \xrightarrow{~\alpha_{d - 1}~} \bigoplus\limits_{k=1}^nB_k \xrightarrow{(b_1,b_2,\cdots,b_n)} X\xrightarrow{~h~} \Sigma^d A}$$
where $b_k\in{\rm rad}_{\C}(B_k,X)$, $k=1,2,\cdots,n$.

Since $\C$ is locally finite, there are only finite many objects $X_i\in\ind(\C),~i=1,2,\cdots,m$
such that $\Hom_{\C}(X_i,X)\neq 0$.
We assume that $\lambda_{ij},~1\leq j\leq q_i$ form a basis of the $k$-vector space ${\rm rad}_{\C}(B_k,X)$.
Put $M:=(\bigoplus\limits_{k=1}^nB_k)\oplus (\bigoplus\limits_{i=1}^n(X_i)^{\oplus q_i})$, we consider the morphism
$$\delta:=(b_1,b_2,\cdots,b_n,\lambda_{11},\cdots,\lambda_{ij},\cdots,\lambda_{mq_m})\in{\rm rad}_{\C}(M,X)$$
which is not retraction, it can be embedded in a $(d+2)$-angle:
$$\xymatrix {M_{\bullet}:~B \xrightarrow{~}M_1 \xrightarrow{~} M_2 \xrightarrow{~} \cdots
  \xrightarrow{~} M_{d-1} \xrightarrow{~} M \xrightarrow{~\delta~} X\xrightarrow{~} \Sigma^d B}.$$
Thus $M_{\bullet}$ is non-split since $\delta$ is not retraction.
Without loss of generality, we can assume that $B=U\oplus V$ where $U$ and $V$ are indecomposable.
By Lemma \ref{y4}, for the morphism $(1,0)\colon U\oplus V\to U$, there exists
the following commutative diagram of $(d+2)$-angles
$$\xymatrix@C=1.2cm{
 U\oplus V \ar[r]^{\;\;(u,~v)}\ar[d]^{(1,~0)}& M_1 \ar[r]\ar@{-->}[d]^{\varphi_1} & M_2 \ar[r]\ar@{-->}[d] & \cdots \ar[r]&M\ar[r]\ar@{-->}[d]&X\ar[r]\ar@{=}[d] & \Sigma^d U\oplus \Sigma^d V \ar[d]^{(1,~0)}\\
U \ar[r]^{\theta_0}&L_1 \ar[r]& L_2 \ar[r] & \cdots \ar[r] & L_d \ar[r]&X\ar[r]^{f}& \Sigma^d U.
}$$
Similarly, for the morphism $(0,1)\colon U\oplus V\to V$, there exists
the following commutative diagram of $(d+2)$-angles
$$\xymatrix@C=1.2cm{
 U\oplus V \ar[r]^{\;\;(u,~v)}\ar[d]^{(0,~1)}& M_1 \ar[r]\ar@{-->}[d]^{\psi_1} & M_2 \ar[r]\ar@{-->}[d] & \cdots \ar[r]&M\ar[r]\ar@{-->}[d]&X\ar[r]\ar@{=}[d] & \Sigma^d U\oplus \Sigma^d V \ar[d]^{(0,~1)}\\
V \ar[r]^{\eta_0}&L_1 \ar[r]& L_2 \ar[r] & \cdots \ar[r] & L_d \ar[r]&X\ar[r]& \Sigma^d V.
}$$
Using similar arguments as in the proof of Lemma \ref{y5}, we conclude
that the at least one of the following two $(d+2)$-angles is non-split
$$\xymatrix{
U \ar[r]^{\theta_0}&L_1 \ar[r] & L_2 \ar[r]  & \cdots \ar[r] & L_d \ar[r]&X\ar[r]^{f\;\;}& \Sigma^d U,}$$
$$\xymatrix{
V \ar[r]^{\eta_0}&N_1 \ar[r] & N_2 \ar[r]  & \cdots \ar[r] & N_d \ar[r]&X\ar[r]& \Sigma^d V.
}$$
Without loss of generality, we assume that
$$\xymatrix{U \ar[r]^{\theta_0}&L_1 \ar[r] & L_2 \ar[r]  & \cdots \ar[r] & L_d \ar[r]&X\ar[r]^{f\;\;\;}& \Sigma^d U,}$$ is non-split.
By Lemma \ref{y3}, we can find a non-split $(d+2)$-angle
$$\xymatrix{
C_{\bullet}:~U \ar[r]^{\quad\omega_0}&C_1 \ar[r]^{\omega_1} & C_2 \ar[r]^{\omega_2}  & \cdots \ar[r]^{\omega_{d-1}} & C_d \ar[r]^{\omega_d}&X\ar[r]^{f\;\;\;}& \Sigma^d U}$$
with $\omega_1,\omega_2,\cdots,\omega_{d-1}$ in ${\rm rad_{\C}}$. Then $C_{\bullet}\in S(X)$.
By (N3), we have the following commutative diagram
$$\xymatrix@C=1.2cm{
 U\oplus V \ar[r]^{\;\;(u,~v)}\ar[d]^{(1,~0)}& M_1 \ar[r]\ar[d]^{\varphi_1} & M_2 \ar[r]\ar[d] & \cdots \ar[r]&M \ar[r]^{\delta}\ar[d]&X\ar[r]\ar@{=}[d] & \Sigma^d U\oplus \Sigma^d V \ar[d]^{(1,~0)}\\
U \ar[r]^{\delta_0}\ar@{=}[d] &L_1 \ar[r]\ar@{-->}[d] & L_2 \ar[r]\ar@{-->}[d]  & \cdots \ar[r] & L_d\ar@{-->}[d] \ar[r]&X\ar[r]^{f}\ar@{=}[d]& \Sigma^d U\ar@{=}[d]\\
U \ar[r]^{\omega_0}&C_1 \ar[r]^{\omega_1} & C_2 \ar[r]^{\omega_2}  & \cdots \ar[r]^{\omega_{d-1}} & C_d \ar[r]^{\omega_d}&X\ar[r]^{f\;\;}& \Sigma^d U
}$$
of $(d+2)$-angles.

For any $D_{\bullet}\in S(X)$, it can be written as
$$\xymatrix {D_{\bullet}:~D\xrightarrow{~}D_1 \xrightarrow{~} D_2 \xrightarrow{~} \cdots
  \xrightarrow{~} D_{d-1} \xrightarrow{~} \bigoplus\limits_{i=1}^pL_i \xrightarrow{d=(d_1,d_2,\cdots,d_p)} X\xrightarrow{~} \Sigma^d D}$$
with $d_i\in{\rm rad}_{\C}(L_i,X)$, $i=1,2,\cdots,p$.
Since $D_{\bullet}\in S(X)$ is non-split, $d$ is not retraction implies that
$d\in{\rm rad}_{\C}(\bigoplus\limits_{i=1}^pL_i ,X)$.
By the definitions of $\lambda_{ij}$ and $\delta$, there exists a morphism
$\rho\colon \bigoplus\limits_{i=1}^pL_i\to M$ such that $d=\delta\rho$.
By (N3), we have the following commutative diagram
$$\xymatrix{
D\ar[r]\ar@{-->}[d]& D_1 \ar[r]\ar@{-->}[d] & D_2 \ar[r]\ar@{-->}[d]& \cdots \ar[r]&\ar[r]\ar@{-->}[d]D_{d-1}&\bigoplus\limits_{i=1}^pL_i \ar[r]^{\quad d}\ar[d]^{\rho}&X\ar[r]\ar@{=}[d] & \Sigma^d D \ar@{-->}[d]\\
B\ar[r]&M_1 \ar[r] & M_2 \ar[r]  & \cdots \ar[r] &M_{d-1}\ar[r]& M \ar[r]^{\delta}&X\ar[r]& \Sigma^d B}$$
of $(d+2)$-angles, where $B=U\oplus V$.
Thus we get the following commutative diagram
$$\xymatrix{
D\ar[r]\ar[d]& D_1 \ar[r]\ar[d] & D_2 \ar[r]\ar[d]& \cdots \ar[r]&\ar[r]\ar[d]D_{d-1}&\bigoplus\limits_{i=1}^pL_i \ar[r]^{\quad d}\ar[d]^{\rho}&X\ar[r]\ar@{=}[d] & \Sigma^d D \ar[d]\\
U \ar[r]^{\omega_0}&C_1 \ar[r]^{\omega_1} & C_2 \ar[r]^{\omega_2}  & \cdots \ar[r]^{\omega_{d-2}}&C_{d-1}\ar[r]^{\omega_{d-1}} & C_d \ar[r]^{\omega_d}&X\ar[r]^{f\;\;}& \Sigma^d U}$$
of $(d+2)$-angles. This shows that $C_{\bullet}$ is a minimal element in $S(X)$.  \qed

\begin{remark}
If there exists a minimal element $S(X)$ or $T(X)$, then it is unique up to isomorphism by Lemma \ref{y3}.
\end{remark}

We are now ready to state and prove our main result.

\begin{theorem}\label{main}
Let $\C$ be a locally finite $(d+2)$-angulated category. If $X\in\ind(\C)$, then there exists an Auslander-Reiten $(d+2)$-angle ending at $X$, and if $X\in\ind(\C)$, then there exists an Auslander-Reiten $(d+2)$-angle starting at $X$.
Thus $\C$ has Auslander-Reiten $(d+2)$-angles.
\end{theorem}

\proof Since $X\in\ind(\C)$, we know that $S(X)$ is non-empty by Lemma \ref{y5}.
By Lemma \ref{y8}, there exists a $(d+2)$-angle
$$A_{\bullet}:~A\xrightarrow{\alpha_0}A_1\xrightarrow{\alpha_1}\cdots\xrightarrow{\alpha_{d-1}}A_{d}\xrightarrow{\alpha_d}X\xrightarrow{\alpha_{d+1}}\Sigma^d A_0$$
which is a minimal element in $S(X)$.
Since $A_{\bullet}\in S(X)$, we have that
$\alpha_1,\alpha_2,\cdots,\alpha_{d_{d-1}}\in {\rm rad}_{\C}$ and $A$ is an indecomposable.
Then $\End(A)$ is local.

We want to prove that $A_{\bullet}$ is an Auslander-Reiten $(d+2)$-angle, by Lemma \ref{lem0},
it suffices to show that $\alpha_d$ is right almost split.

Assume that $g\colon M\to X$ is not retraction.
By the dual of Lemma \ref{y4}, there exists
the following commutative diagram of $(d+2)$-angles
$$\xymatrix@C=1.2cm{
A\ar[r]^{\gamma_0}\ar@{=}[d]& B_1 \ar[r]^{\gamma_1}\ar@{-->}[d]^{\psi_1} & B_2 \ar[r]^{\gamma_2}\ar@{-->}[d]^{\psi_2}  & \cdots \ar[r]^{\gamma_{d-1}}& B_d \ar[r]^{\gamma_d\;\;}\ar@{-->}[d]^{\psi_d}&M\ar[r]^{\gamma_{d+1}}\ar[d]^{\gamma_d} & \Sigma^d A \ar@{=}[d]\\
A\ar[r]^{\alpha_0}&A_1 \ar[r]^{\alpha_1} & A_2 \ar[r]^{\alpha_2}  & \cdots \ar[r]^{\alpha_{d-1}} & A_d \ar[r]^{\alpha_d}&X\ar[r]^{\alpha_{d+1}}& \Sigma^d A
}$$
such that
$$N_{\bullet}:~B_1\xrightarrow{} N_1\xrightarrow{}
 N_2\xrightarrow{}\cdots\xrightarrow{}N_{d-1}
                                     \xrightarrow{}M\oplus A_d\xrightarrow{\left(
                                       \begin{smallmatrix}
                                         g,&\alpha_{d}
                                       \end{smallmatrix}
                                     \right)}X\xrightarrow{~h~}\Sigma^dB_1$$
 is a $(d+2)$-angle in $\C$, where $N_i=B_{i+1}\oplus A_{i},~i=1,2,\cdots,d-1$.
Since $g$ and $\alpha_d$ are not retraction, we have that
$(g,\alpha_d)$ is also not retraction by using similar arguments as in the proof of Lemma \ref{y7}.
That is, $N_{\bullet}$ is non-split.

Without loss of generality, we can assume that $B_1=U\oplus V$ where $U$ and $V$ are indecomposable.
By Lemma \ref{y4}, for the morphism $(1,0)\colon U\oplus V\to U$, there exists
the following commutative diagram of $(d+2)$-angles
$$\xymatrix@C=1.2cm{
 U\oplus V \ar[r]^{\;\;(u,~v)}\ar[d]^{(1,~0)}& N_1 \ar[r]\ar@{-->}[d]^{\varphi_1} & N_2 \ar[r]\ar@{-->}[d] & \cdots \ar[r]&M\oplus A_d \ar[r]\ar@{-->}[d]&X\ar[r]\ar@{=}[d] & \Sigma^d U\oplus \Sigma^d V \ar[d]^{(1,~0)}\\
U \ar[r]^{\delta_0}&L_1 \ar[r]& L_2 \ar[r] & \cdots \ar[r] & L_d \ar[r]&X\ar[r]^{f\;\;}& \Sigma^d U.
}$$
Similarly, for the morphism $(0,1)\colon U\oplus V\to V$, there exists
the following commutative diagram of $(d+2)$-angles
$$\xymatrix@C=1.2cm{
 U\oplus V \ar[r]^{\;\;(u,~v)}\ar[d]^{(0,~1)}& N_1 \ar[r]\ar@{-->}[d]^{\psi_1} & N_2 \ar[r]\ar@{-->}[d] & \cdots \ar[r]&M\oplus A_d \ar[r]\ar@{-->}[d]&X\ar[r]\ar@{=}[d] & \Sigma^d U\oplus \Sigma^d V \ar[d]^{(0,~1)}\\
V \ar[r]^{\eta_0}&Q_1 \ar[r]& Q_2 \ar[r] & \cdots \ar[r] & Q_d \ar[r]&X\ar[r]& \Sigma^d V.
}$$
Using similar arguments as in the proof of Lemma \ref{y5}, we conclude
that the at least one of the following two $(d+2)$-angles is non-split
$$\xymatrix{
U \ar[r]^{\delta_0}&L_1 \ar[r] & L_2 \ar[r]  & \cdots \ar[r] & L_d \ar[r]&X\ar[r]^{f\;\;}& \Sigma^d U,}$$
$$\xymatrix{
V \ar[r]^{\eta_0}&Q_1 \ar[r] & Q_2 \ar[r]  & \cdots \ar[r] & Q_d \ar[r]&X\ar[r]& \Sigma^d V.
}$$
Without loss of generality, we assume that
$$\xymatrix{U \ar[r]^{\delta_0}&L_1 \ar[r] & L_2 \ar[r]  & \cdots \ar[r] & L_d \ar[r]&X\ar[r]^{f\;\;}& \Sigma^d U,}$$ is non-split.
By Lemma \ref{y3}, we can find a non-split $(d+2)$-angle
$$\xymatrix{
C_{\bullet}:~U \ar[r]^{\quad\lambda_0}&C_1 \ar[r]^{\lambda_1} & C_2 \ar[r]^{\lambda_2}  & \cdots \ar[r]^{\lambda_{d-1}} & C_d \ar[r]^{\lambda_d}&X\ar[r]^{f\;\;}& \Sigma^d U}$$
with $\lambda_1,\lambda_2,\cdots,\lambda_{d-1}$ in ${\rm rad_{\C}}$.
By (N3), we have the following commutative diagram
$$\xymatrix@C=1.2cm{A\ar[r]^{\alpha_0}\ar@{-->}[d]& A_1 \ar[r]^{\alpha_1}\ar@{-->}[d] & A_2 \ar[r]^{\alpha_2}\ar@{-->}[d]& \cdots \ar[r]^{\alpha_{d-1}}& A_d \ar[r]^{\alpha_d\;\;}\ar[d]^{\binom{0}{1}}&X\ar[r]^{\alpha_{d+1}}\ar@{=}[d] & \Sigma^d A_0 \ar@{-->}[d]\\
 U\oplus V \ar[r]^{\;\;(u,~v)}\ar[d]^{(1,~0)}& N_1 \ar[r]\ar[d]^{\varphi_1} & N_2\ar[r]\ar[d] & \cdots \ar[r]&M\oplus A_d \ar[r]^{\quad(g,\alpha_d)}\ar[d]&X\ar[r]\ar@{=}[d] & \Sigma^d U\oplus \Sigma^d V \ar[d]^{(1,~0)}\\
U \ar[r]^{\delta_0}\ar@{=}[d] &L_1 \ar[r]\ar@{-->}[d] & L_2 \ar[r]\ar@{-->}[d]  & \cdots \ar[r] & L_d\ar@{-->}[d] \ar[r]&X\ar[r]^{f}\ar@{=}[d]& \Sigma^d U\ar@{=}[d]\\
U \ar[r]^{\lambda_0}&C_1 \ar[r]^{\lambda_1} & C_2 \ar[r]^{\lambda_2}  & \cdots \ar[r]^{\lambda_{d-1}} & C_d \ar[r]^{\lambda_d}&X\ar[r]^{f\;\;}& \Sigma^d U
}$$
of $(d+2)$-angles.
We obtain that $A_{\bullet}>C_{\bullet}$ implies that $A_{\bullet}\sim C_{\bullet}$ since $A_{\bullet}$ is the minimal element in $S(X)$. Thus there exists the following commutative diagram
$$\xymatrix@C=1.2cm{U \ar[r]^{\lambda_0}\ar[d]&C_1 \ar[r]^{\lambda_1} \ar[d]& C_2 \ar[r]^{\lambda_2}\ar[d]  & \cdots \ar[r]^{\lambda_{d-1}} & C_d\ar[d] \ar[r]^{\lambda_d}&X\ar[r]^{f\;\;}\ar@{=}[d]& \Sigma^d U\ar[d]\\
A\ar[r]^{\alpha_0}& A_1 \ar[r]^{\alpha_1} & A_2 \ar[r]^{\alpha_2}& \cdots \ar[r]^{\alpha_{d-1}}& A_d \ar[r]^{\alpha_d\;\;}&X\ar[r]^{\alpha_{d+1}}& \Sigma^d A_0
}$$
of $(d+2)$-angles.
Hence we get the following commutative diagram
$$\xymatrix@C=1.2cm{ U\oplus V \ar[r]^{\;\;(u,~v)}\ar[d]& N_1 \ar[r]\ar[d] & N_2\ar[r]\ar[d] & \cdots \ar[r]&M\oplus A_d \ar[r]^{\quad(g,\alpha_d)}\ar[d]^{(a,b)}&X\ar[r]\ar@{=}[d] & \Sigma^d U\oplus \Sigma^d V \ar[d]\\
A\ar[r]^{\alpha_0}& A_1 \ar[r]^{\alpha_1} & A_2 \ar[r]^{\alpha_2}& \cdots \ar[r]^{\alpha_{d-1}}& A_d \ar[r]^{\alpha_d\;\;}&X\ar[r]^{\alpha_{d+1}}& \Sigma^d A_0
}$$
of $(d+2)$-angles. It follows that $g=a\alpha_d$.
This shows that $\alpha_d$ is right almost split.

Similarly, we can show that if $X\in\ind(\C)$, then there exists an Auslander-Reiten $(d+2)$-angle starting at $X$.
Thus $\C$ has Auslander-Reiten $(d+2)$-angles.  \qed

\begin{remark}
As a special case of Theorem \ref{main} when $d=1$, that is, if $\C$ is a locally finite
triangulated category, then $\C$ has Auslander-Reiten triangles, see \cite{xz1,xz2}.
\end{remark}

\textbf{Panyue Zhou}\\
College of Mathematics, Hunan Institute of Science and Technology, 414006, Yueyang, Hunan, P. R. China.\\
and \\
D\'{e}partement de Math\'{e}matiques, Universit\'{e} de Sherbrooke, Sherbrooke,
Qu\'{e}bec J1K 2R1, Canada.\\
E-mail: \textsf{panyuezhou@163.com}

\end{document}